\documentclass[12pt]{article}
\usepackage{amssymb}

\newtheorem{thm}     {Theorem}[section]

\newtheorem{definition}  [thm]{Definition}

\newtheorem{lemma}   [thm]{Lemma}

\newcommand{\proof} {\noindent{\bf Proof. }}

\newcommand{\B}{\mathbb B}
\newcommand{\C}{\mathbb C}
\newcommand{\D}{\mathbb D}

\newcommand{\R}{\mathbb R}

\newcommand{\st}{{\rm st}}

\def\Im{{\rm Im\,}}
\def\bar{\overline}

\begin{document}

\title{The Chirka - Lindel\"of and Fatou theorems for $\overline\partial_J$-subsolutions }
\author{Alexandre Sukhov{*} }
\date{}
\maketitle

{\small

* Universit\'e  de Lille, Laboratoire
Paul Painlev\'e,
U.F.R. de
Math\'e-matique, 59655 Villeneuve d'Ascq, Cedex, France, sukhov@math.univ-lille1.fr
The author is partially suported by Labex CEMPI.

Institut of Mathematics with Computing Centre - Subdivision of the Ufa Research Centre of Russian
Academy of Sciences, 45008, Chernyshevsky Str. 112, Ufa, Russia.

}
\bigskip

{\small Abstract.  This paper studies boundary properties of bounded functions with bounded 
$\overline\partial_J$ differential on strictly pseudoconvex domains in an almost complex manifold.}



\bigskip

\tableofcontents

\section{Introduction}
The first fundamental results on analytic properties of almost complex structures (in several variables) are due to Newlander - Nirenberg \cite{NeNi} and Nijenhuis - Woolf \cite{NiWo}. After the seminal work by M.Gromov \cite{Gr} the theory of pseudoholomorphic curves in almost complex manifolds became one of the most powerful tools of the symplectic geometry and  now is rapidly increasing.

According to the Nijenhuis - Woolf theorem, every almost complex manifold $(M,J)$ locally (near every point) contains  plenty of pseudoholomorphic curves i.e. open Riemann surfaces embedded to $M$ compatibly with an almost complex structure $J$. Another object which always exists locally is (strictly) plurisubharmonic functions; they  in turn are indespensable for symplectic applications (also being of independent interest, of course). For these reasons the study of pseudoholomorphic curves and plurisubharmonic functions have many common points with the case of usual complex manifolds (i.e. with the case of integrable complex structures), although the almost complex proofs often require much more involved tools from elliptic PDEs and the non-linear analysis. A serious difference with respect to the integrable case arises when one tries to consider holomorphic functions on an almost complex manifold i.e. the solutions $f$ of the equation $\overline\partial_J f = 0$. In the case of complex dimension $ > 1$, for a "generically choosen" almost complex structure $J$ all (even only locally defined!) holomorphic functions are constant. Thus an attempt to  extend directly the function theory in several variables to the almost complex case is fruitless.

A remarkable progress here was done by S.Donaldson \cite{D}. Here considered the almost holomorphic functions $f$ which locally satisfy the Beltrami type condition $\parallel \overline\partial_J f \parallel \le \varepsilon \parallel \partial_J f \parallel$ with $\varepsilon < 1$ (one can say that these functions are subsolutions of the multidimensional Beltrami operator). Such objects always exist locally for each almost complex structure $J$. If $J$ is compatible with some symplectic structure, the zero set of an almost holomorphic function is a symplectic hypersurface. Using these tools, S. Donaldson obtain a symplectic analog of the Lefschetz hyperplane section theorem having fundamental consequences in the symplectic geometry. The approach of Donaldson provides the second Complex Analysis tool for the symplectic geometry, quite independently from Gromov's theory. This is a partial motivation of the present paper.

Donaldson's work opens a natural way to develop a function theory on almost complex manifolds. The idea is to study classes of functions imposing suitable assumptions on the  $\overline\partial_J$-part of differential . The assumptions depend on the problem under investigation. For example, local properties of subsolutions to the Beltrami operator have some common points with usual holomorphic functions and are in the focus of Donaldson's theory. In the present paper we study boundary properties and choose other type of conditions. This is well-known in the function theory that many boundary properties of functions $f$  in domains of $\C^n$ with some boundary control over the $\overline\partial f$ (i.e. the subsolutions to the $\overline\partial$-equation) are similar to the ones of usual holomorphic functions. Sometimes such functions are called asymptotically holomorphic; they were succesfully applied, for instance,  in the works of \cite{NiWeY, NaRu, Ro, PiKh,Fo} and in many others. The goal of the present paper is to study some boundary properties of their almost complex analogs.  It is worth to stress that locally there are plenty of such functions on any almost complex manifold and so they  represent a natural object for study.

There exist several ways to define the boundary control over the $\overline\partial_J f$ and their choice usually depend on boundary properties which are under the study. In the present paper  we establish the Chirka - Lindel\"of and Fatou type theorems for bounded $C^1$-functions $f$ with bounded $\overline\partial_J f$ on a strictly pseudoconvex domain in an almost complex manifold $(M,J)$ (see resp. Theorem \ref{Thm1} and \ref{Thm2}). The main results are contained in Sections 3 and 4; we potspone a detailed duscussion of our methods and relations with previous results until these sections in order to avoid a long introduction. Here I only mention that the main results are inspired by the works and methods of E.Stein \cite{St}, E.Chirka \cite{Ch} and A.Sadullaev \cite{Sa1} which now became classical; one can consider the obtained results as their  generalization.

\section{ Almost complex manifolds  and almost holomorphic functions} In this rather long preliminary section we recall basic notions of the almost complex geometry making our presentation more convenient for specialists in Analysis. A reader could find much more detailed information in \cite{Aud}.
Everywhere through this paper we assume that manifolds and almost complex structures are of class $C^\infty$ (the word "smooth" means the regularity of this class); notice the main results remain true under considerably weaker regularity assumptions.

\subsection{ Almost complex manifolds} Let $M$ be a smooth  manifold of dimension $2n$. {\it An almost complex structure} $J$ on $M$ is a smooth map  which associates to every point $p \in M$ a linear isomorphism $J(p): T_pM \to T_pM$ of the tangent space $T_pM$ such that $J(p)^2 = -I$; here  $I$ denotes the identity map of $T_pM$. Thus, every linear map $J(p)$ is a  complex structure on a vector space $T_pM$ in the usual sense of Linear Algebra. A couple $(M,J)$ is called {\it an almost complex manifold} of complex dimension n. Note that every almost complex manifold admits the canonical orientation represented by $(e_1,Je_1,....,e_n,Je_n)$ where $(e_1,....,e_n)$ is any complex basis of $(T_pM,J(p))$.

One of the most  important examples is provided by the {\it standard complex structure} $J_{st} = J_{st}^{(2)}$ on $M = \R^2$; it is represented in the canonical coordinates of $\R^2$ by the matrix 

\begin{eqnarray}
\label{J_st}
J_{st}= \left(
\begin{array}{cll}
0 & & -1\\
1 & & 0
\end{array}
\right)
\end{eqnarray}
More generally, the standard complex structure $J_{st}$ on $\R^{2n}$ is represented by the block diagonal matrix $diag(J_{st}^{(2)},...,J_{st}^{(2)})$ (usually we drop the notation of dimension because its value  will be clear from the context). As usual,  setting $iv := Jv$ for $v \in \R^{2n}$, we identify $(\R^{2n},J_{st})$ with $\C^n$; we  use the notation 
$z = x + iy = x + Jy$ for the standard complex coordinates $z = (z_1,...,z_n) \in \C^n$.

Let $(M,J)$ and $(M',J')$ be smooth  almost complex manifolds. A $C^1$-map $f:M' \to M$ is called  
$(J',J)$-complex or  $(J',J)$-holomorphic  if it satisfies {\it the Cauchy-Riemann equations} 
\begin{eqnarray}
\label{CRglobal}
df \circ J' = J \circ df.
\end{eqnarray}

This is easy to check that   a map $f: \C^n \to \C^m$ is $(J_{st},J_{st})$-holomorphic if and only if each component of $f$ is a usual holomorphic function.

 Every almost complex manifold
$(M,J)$ can be viewed locally as the unit ball $\B$ in
$\C^n$ equipped with a small (in any $C^m$-norm) almost complex
deformation of $J_{st}$. The following statement is often very useful.
\begin{lemma}
\label{lemma1}
Let $(M,J)$ be an almost complex manifold. Then for every point $p \in
M$, every  $m \geq 0$ and   $\lambda_0 > 0$ there exist a neighborhood $U$ of $p$ and a
coordinate diffeomorphism $z: U \rightarrow \B$ such that
$z(p) = 0$, $dz(p) \circ J(p) \circ dz^{-1}(0) = J_{st}$,  and the
direct image $ z_*(J): = dz \circ J \circ dz^{-1}$ satisfies $\vert\vert z_*(J) - J_{st}
\vert\vert_{C^m(\bar {\B})} \leq \lambda_0$.
\end{lemma}
\proof There exists a diffeomorphism $z$ from a neighborhood $U'$ of
$p \in M$ onto $\B$ satisfying $z(p) = 0$; after an additional linear change of coordinates one can achieve  $dz(p) \circ J(p)
\circ dz^{-1}(0) = J_{st}$ (this is a classical fact from the Linear Algebra). For $\lambda > 0$ consider the isotropic dilation
$h_{\lambda}: t \mapsto \lambda^{-1}t$ in $\R^{2n}$ and the composition
$z_{\lambda} = h_{\lambda} \circ z$. Then $\lim_{\lambda \rightarrow
0} \vert\vert (z_{\lambda})_{*}(J) - J_{st} \vert\vert_{C^m(\bar
{\B})} = 0$ for every  $m \geq 0$. Setting $U = z^{-1}_{\lambda}(\B)$ for
$\lambda > 0$ small enough, we obtain the desired statement. 
In what follows we often denote the structure $z^*(J)$ again by $J$ viewing it as a local representation of $J$ in the coordinate system $(z)$.

Recall that an almost complex structure $J$ is called {\it integrable} if $(M,J)$ is locally biholomorphic in a neighborhood of each point to an open subset of $(\C^n,J_{st})$. In the case of complex dimension 1 every almost complex structure is integrable. In the case of complex dimension $> 1$  integrable almost complex structures form a highly special subclass in the space of all almost complex structures on $M$; an efficient criterion of integrablity is provided by the classical theorem of Newlander - Nirenberg \cite{NeNi}.

\bigskip

 \subsection{Pseudoholomorphic discs} Let $(M,J)$ be an almost complex manifold of dimension $n > 1$. For a "generic" choice of an almost complex structure, any holomorphic (even locally) function on $M$ is constant. Similarly, $M$ does not admit non-trivial $J$-complex submanifolds (that is, submanifolds with tangent spaces invariant with respect to  $J$) of complex dimension $> 1$. The only (but fundamentally important) exception arises in the case of pseudoholomorphic curves i.e. $J$-complex submanifolds of  complex dimension 1: they always exist locally.

Usually  pseudoholomorphic curves arise in connection with   solutions  $f$ of (\ref{CRglobal}) in the special case   where $M'$ has the complex dimension 1. These holomorphic maps are called $J$-complex (or $J$-holomorphic or {\it pseudoholomorphic} ) curves. Note that we view here the curves as maps i.e. we consider parametrized curves.
We use the notation  $\D = \{ \zeta \in \C: \vert \zeta \vert < 1 \}$ for  the
unit disc in $\C$ always assuming that it is equipped with the standard complex structure   $J_{\st}$. If in the equations (\ref{CRglobal})  we have $M' = \D$ ; we  call such a map $f$ a $J$-{\it complex  disc} or a  {\it pseudoholomorphic disc} or just a  holomorphic disc
when the structure  $J$ is fixed. 

A fundamental fact is that  pseudoholomorphic discs always exist in a suitable neighborhood of any point of $M$; this is the classical Nijenhuis-Woolf theorem (see \cite{NiWo}). Here it is convenient to rewrite the equations (\ref{CRglobal}) in local coordinates  similarly to the complex version of the usual Cauchy-Riemann equations.

 Everything will be local, so (as above) we are in a neighborhood $\Omega$ of $0$ in $\C^n$ with the standard complex coordinates $z = (z_1,...,z_n)$. We assume that $J$ is an almost complex structure defined on $\Omega$ and $J(0) = J_{st}$. Let 
$$z:\D \to \Omega,$$ 
$$z : \zeta \mapsto z(\zeta)$$ 
be a $J$-complex disc. Setting $\zeta = \xi + i\eta$ we write (\ref{CRglobal}) in the form $z_\eta = J(Z) Z_\xi$. This equation can be in turn written as

\begin{eqnarray}
\label{holomorphy}
z_{\bar\zeta} - A(z)\bar z_{\bar\zeta} = 0,\quad
\zeta\in\D.
\end{eqnarray}
Here a smooth map $A: \Omega \to Mat(n,\C)$ is defined by the equality $L(z) v = A \overline v$ for any vector $v \in \C^n$ and $L$ is an $\R$-linear map defined by $L = (J_{st} + J)^{-1}(J_{st} - J)$. It is easy to check that the condition $J^2 = -Id$ is equivalent to the fact that $L$ is $\overline\C$-linear. The matrix $A(z)$ is called {\it the complex matrix} of $J$ in the local coordinates $z$. Locally the correspondence between $A$ and $J$ is one-to-one. Note that the condition $J(0) = J_{st}$ means that $A(0) = 0$. 

If $z'$ are other local coordinates and $A'$ is the corresponding complex matrix of $J'$, then, as it is easy to check, we have the following transformation rule:

\begin{eqnarray}
\label{CompMat}
A' = (z_z' A  + { z}_{\overline z}')({\overline z}_{\overline z}' + {\overline z}_{ z}'A)^{-1}
\end{eqnarray}
(see \cite{SuTu}).

Note that one can view the equations (\ref{holomorphy}) as a quasilinear analog of the Beltrami equation for vector-functions. From this point of view, the theory of pseudoholomorphic curves is an analog of the theory of quasi-conformal mappings.

\bigskip

Recall that for a complex function $f$  {\it the Cauchy-Green transform} is defined by

\begin{eqnarray}
\label{CauchyGreen}
Tf(\zeta) = \frac{1}{2 \pi i} \int\int_{\D} \frac{f(\omega)d\omega \wedge d\overline\omega}{\omega - \zeta}
\end{eqnarray}
This is the main analytic tool in the theory of pseudoholomorphic curves. This is classical that the operator $T$ has the following properties:
\begin{itemize}
\item[(i)] $T: C^r(\D) \to C^{r+1}(\D)$ is a bounded linear operator for every non-integer $r > 0$ ( a similar property holds in the Sobolev scales, see below). Here we use the usual H\"older norm on the space $C^r(\D)$.
\item[(ii)] $(Tf)_{\overline\zeta} = f$ i.e. $T$ solves the $\overline\partial$-equation in the unit disc. 
\item[(iii)] the function $Tf$ is holomorphic on $\C \setminus \overline\D$.
\end{itemize}
Fix a real non-integer $r > 1$. Let $z: \D \to \C^n$, $z: \D \ni \zeta \mapsto z(\zeta)$ be a $J$-complex disc. 
Since  the operator
$$\Psi_{J}: z \longrightarrow w =  z - TA(z) \overline {z}_{\overline \zeta} $$
takes the space   $C^{r}(\overline{\mathbb D})$  into itself,  we can write   the
equation (\ref{CRglobal}) in the form 
$(\Psi_J(z))_{\overline \zeta}  = 0$. Thus, the disc $z$ is $J$-holomorphic if
and only if the map $\Psi_{J}(z):\mathbb D \longrightarrow \C^n$ is
$J_{st}$-holomorphic.
When the norm of $A$  is small enough (which is assured  by Lemma \ref{lemma1}),
then  by the implicit function theorem the operator    $\Psi_J$
is invertible  and we obtain a bijective
correspondence between  $J$-holomorphic discs and usual
holomorphic discs. This easily implies the existence of a $J$-holomorphic disc
in a given tangent direction through a given point of $M$, as well as  a smooth dependence of such a
disc  on a deformation of a point or a tangent vector, or on an almost complex structure; this also establishes  the interior elliptic regularity of discs.

\bigskip

Let $(M,J)$ be an almost complex manifold and $E \subset M$ be a real submanifold of $M$. 
Suppose that a $J$-complex disc $f:\D \to M$ is  continuous on $\overline\D$.  We some abuse of terminology, we also call the image $f(\D)$  simply by a disc and twe call he image $f(b\D)$ by  the boundary of a disc. If  $f(b\D) \subset E$, then we say that (the boundary of ) the disc  $f$ is {\it glued} or {\it attached} to $E$ or simply 
that $f$ is attached to $E$. If $\Gamma \subset b\D$ is an arc and $f(\Gamma) \subset E$, we say that $f$ is glued or attached to $E$ along $\Gamma$.

\subsection{The $\overline\partial_J$-operator on an almost complex manifold $(M,J)$}

Consider now the second special class (together with pseudoholomorphic curves) of holomorphic maps. Consider first the situation when $J$ be an almost complex structure defined 
in a domain $\Omega\subset\C^n$; one can view this as a local coordinate representation of $J$ in a chart on $M$.

A $C^1$ function $F:\Omega\to\C$ is $(J,J_{st})$-holomorphic
if and only if it satisfies the Cauchy-Riemann equations
\begin{eqnarray}
\label{CRscalar}
F_{\bar z} + F_z A(z)  =0,
\end{eqnarray}
where $F_{\bar z} = (\partial F/\partial \overline{z}_1,...,\partial F/\partial \overline{z}_n)$ and $F_z = (\partial F/\partial {z}_1,...,\partial F/\partial {z}_n)$ are viewed as  row-vectors. Indeed, $F$ is $(J,J_{st})$ holomorphic if and only if 
for every $J$-holomorphic disc $z:\D \to \Omega$ the composition $F \circ z$ is a usual holomorphic function  that is $\partial (F \circ z) /\partial\overline{\zeta} = 0$ on $\D$. Then the  Chain rule in combination with (\ref{holomorphy}) leads to (\ref{CRscalar}). Generally the only solutions to (\ref{CRscalar}) are constant functions  unless $J$ is integrable (then $A$ vanishes identically in suitable coordinates). Note also that (\ref{CRscalar}) is a linear PDE system while (\ref{holomorphy}) is a quasilinear PDE for a vector function on $\D$.

Every $1$-differnitial form $\phi$ on $(M,J)$ admits a unique decomposition $\phi = \phi^{1,0} + \phi^{0,1}$ with respect to $J$. In particular, if $F:(M,J) \to \C$ is a $C^1$-complex function, we have $dF = dF^{1,0} + dF^{0,1}$. We use the notation 
\begin{eqnarray}
\label{d-bar}
\partial_J F = dF^{1,0} \,\,\,\mbox{and}\,\,\, \overline\partial_J F = dF^{0,1}
\end{eqnarray}

In order to write these operators explicitely in local coordinates, we find a  local basic in the space of (1,0) and (0,1) forms. We view  $dz = (dz_1,...,dz_n)^t$ and $d\overline{z} = (d\overline{z}_1,...,d\overline{z}_n)^t$ as vector-columns. Then the forms 
\begin{eqnarray}
\label{FormBasis}
\alpha = (\alpha_1,..., \alpha_n)^t = dz - A d\overline{z} \,\,\, \mbox{and} \,\, \overline\alpha = d\overline{z} - \overline A dz
\end{eqnarray}
form a basis in the space of  (1,0) and (0,1) forms respectively. Indeed, it suffices to note that for 1-form $\beta$ is (1,0) (resp. $(0,1)$) for if and only if for every $J$-holomorphic disc $z$ the pull-back $z^*\beta$ is a usual (1,0) (resp. $(0,1)$) form on $\D$. Using the equations (\ref{holomorphy}) we obtain the claim.

Now we decompose the differential $dF = F_zdz + F_{\overline{z}} d\overline{z} = \partial_J F + \overline\partial_J F$ in the basis $\alpha$, $\overline\alpha$ using (\ref{FormBasis}) and obtain the explicit expression 
\begin{eqnarray}
\label{d-bar2}
 \overline\partial_J F = (F_{\overline{z}} (I - \overline{A}A)^{-1} + F_z (I - A\overline{A})^{-1}A)\overline\alpha
\end{eqnarray}

It is easy to check that the holomorphy condition $\overline\partial_J F = 0$ is equivalent to (\ref{CRscalar}) because $(I - A\overline{A})^{-1} A (I - \overline{A} A) = A$. Thus 

\begin{eqnarray*}
 \overline\partial_J F = (F_{\overline{z}}  + F_z A)(I - \overline{A}A)^{-1}\overline\alpha
\end{eqnarray*}

We note that the term $(I - A\overline A)^{-1}$ as well as the forms $\alpha$ affect only the non-essential constants in local estimates of the $\overline\partial_J$-operator near a boundary point which we will perfom in the next sections. So the reader  can assume that this operator is simply given by the left hand of (\ref{CRscalar}).

Let $F$ be a complex function of class $C^1$ on a (bounded) domain $\Omega$ in an almost complex manifold of dimension $n$. We call a function $f$ a {\it subsolution of the $\overline\partial_J$ operator }  or simply 
$\overline\partial_J$-subsolution  on $\Omega$ if $\parallel \overline\partial_J F \parallel $ is uniformly bounded on $\Omega$ that is there exists a constant $C > 0$ such that
\begin{eqnarray}
\label{AlHol}
\parallel \overline\partial_J F (z)\parallel \le C
\end{eqnarray}
for all $z \in \Omega$. Here we use the norm with respect to any fixed Riemannian metric on $M$.

Obviously, non-constant $\overline\partial_J$-subsolutions exist in a sufficiently small neighborhhod of any point of $M$. In fact any function $F$ of class $C^1$ in an open neighborhhod of the compact set $\overline\Omega$ is a $\overline\partial_J$-subsolution on $\Omega$.

Let $F$ be a $\overline\partial_J$-subsolution on $\Omega$. Suppose that $A$ is the complex matrix of $J$ in a local chart $U$ and $z:\D \to U$ is a $J$-complex disc. It follows by the Chain Rule and (\ref{holomorphy}) that
$$(F \circ z)_{\overline\zeta} = (F_{\overline z} + F_zA){\overline z}_{\overline\zeta}.$$
Thus, if $h: \D \to \Omega$ is a $J$-complex disc of class $C^1(\overline\D)$, then the composition 
$F \circ h$ has a uniformly bounded $\overline\partial$-derivative on $\D$ that is $F \circ h$ is a $\overline\partial_{J_{st}}$-subsolution on $\D$. Note that the upper bound on the $\overline\partial (F \circ h)$ depends only on the upper bound on $\overline\partial_J F$ from (\ref{AlHol}) and the $C^1$ norm of $h$ on $\overline\D$. In particular, if $(h_t)$ is a family of $J$-complex discs in $\Omega$ and $C^1$-norms of these discs are uniformly bounded  with respect $t$, then the norms $\parallel\overline\partial_J (F \circ h_t) \parallel$ are bounded uniformly in $t$ as well.

\subsection{Plurisubharmonic functions on almost complex manifolds: the background} For the convenience of readers we recall some basic notions concerning plurisubharmonic functions on almost complex manifolds.  Let  $u$ be a real $C^2$ function on an open subset $\Omega$ of an almost complex manifold  $(M,J)$. Denote by $J^*du$ the
 differential form acting on a vector field $X$ by $J^*du(X):= du(JX)$. Given point $p \in M$ and a tangent vector $V \in T_p(M)$ consider  a smooth vector field $X$ in a
neighborhood of $p$ satisfying $X(p) = V$. 
The value of the {\it complex Hessian} ( or  the  Levi form )   of $u$ with respect to $J$ at $p$ and $V$ is defined by $H(u)(p,V):= -(dJ^* du)_p(X,JX)$.  This definition is independent of
the choice of a vector field $X$. For instance, if $J = J_{st}$ in $\C$, then
$-dJ^*du = \Delta u d\xi \wedge d\eta$; here $\Delta$ denotes the Laplacian. In
particular, $H_{J_{st}}(u)(0,\frac{\partial}{\partial \xi}) = \Delta u(0)$.

Recall some  basic properties of the complex Hessian (see for instance, \cite{DiSu}):

\begin{lemma}
\label{pro1}
Consider  a real function $u$  of class $C^2$ in a neighborhood of a point $p \in M$.
\begin{itemize}
\item[(i)] Let $F: (M',J') \longrightarrow (M, J)$ be a $(J',
  J)$-holomorphic map, $F(p') = p$. For each vector $V' \in T_{p'}(M')$ we have
$H_{J'}(u \circ F)(p',V') = H_{J}(u)(p,dF(p)(V'))$.
\item[(ii)] If $f:\D \longrightarrow M$ is a $J$-complex disc satisfying
  $f(0) = p$, and $df(0)(\frac{\partial}{\partial \xi}) = V \in T_p(M)$ , then $H_J(u)(p,V) = \Delta (u \circ f) (0)$.
\end{itemize}
\end{lemma}
 Property (i) expresses the holomorphic invariance of the complex Hessian. Property (ii) is often useful in order to compute the complex Hessian  on a given 
tangent vector $V$.

 Let $\Omega$ be a domain $M$. An upper semicontinuous function $u: \Omega \to [-\infty,+\infty[$ on $(M,J)$ is
{\it $J$-plurisubharmonic} (psh) if for every $J$-complex disc $f:\D \to \Omega$ the composition $u \circ f$ is a subharmonic function on $\D$. By Proposition
\ref{pro1}, a $C^2$ function $u$ is psh on $\Omega$ if and only if it has    a 
positive semi-definite complex Hessian on $\Omega$ i.e. $H_J(u)(p,V) \geq 0$  for any $ p \in \Omega $ and $V \in T_p(M)$. 
A real $C^2$ function $u:\Omega \to \R$ is called {\it strictly $J$-psh} on $\Omega$, if $H_J(u)(p,V) > 0$ for each $p \in M$ and $V \in T_p(M) \backslash \{ 0\}$. Obviously, these notions   are local: an upper semicontinuous (resp. of class $C^2$) function on $\Omega$ is  $J$-psh (resp. strictly) on $\Omega$ if and only if it is $J$-psh (resp. strictly) in some open neighborhood of each point of $\Omega$.

A useful  observation is that the Levi form of a function $r$ at
a point $p$ in
an almost complex manifold $(M,J)$ coincides with the Levi form with respect
to the standard structure $J_{st}$ of $\R^{2n}$ if {\it suitable} local
coordinates near $p$ are choosen. Let us explain how to construct these adapted
coordinate systems.

As above, choosing local coordinates near $p$ we may identify a neighborhood
of $p$ with a neighborhood of the origin and assume that $J$-holomorphic discs
are solutions of (\ref{holomorphy}).

\begin{lemma}
\label{normalization}
There exists a  second order polynomial local diffeomorphism $\Phi$ fixing the
origin and with linear part equal to the identity such that in the new coordinates
 the complex matrix  $A$  of $J$ (that is $A$ from the equation (\ref{holomorphy})) satisfies
\begin{eqnarray}
\label{norm}
A(0) = 0, A_{z}(0) = 0
\end{eqnarray}
\end{lemma}
Thus, by a suitable local change of coordinates one can remove the 
terms linear in $z$ in the matrix $A$. We stress that in general it is impossible to
get rid of  first order terms containing  $\overline z$ since this would
impose a restriction on the Nijenhuis tensor $J$ at the origin.

I have learned this result  from  unpublished E.Chirka's notes; see .\cite{DiSu} for the proof. In
\cite{SuTu} it is shown that, in an almost complex manifold of (complex) dimension 2,
 a similar normalization is possible along a given embedded $J$-holomorphic disc.

\subsection{Boundary properties of subsolutions of the $\overline\partial$-operator  in the unit disc}

Denote by $W^{k,p}(\D)$ (we need only $k=0$ and $k=1$) the usual Sobolev classes of functions having generalized partial derivatives up to the order $k$ in $L^p(\D)$ ; thus, $W^{0,p}(\D) = L^p(\D)$. We will always assume that $p > 2$.

Consider {\it the Cauchy transform}

\begin{eqnarray}
\label{Cauchy}
Kf(\zeta) = \frac{1}{2 \pi i} \int_{b\D} \frac{f(\omega)d\omega \wedge d\overline\omega}{\omega - \zeta}
\end{eqnarray}
Recall that $K$ is a bounded linear map in classes  $C^r(b\D) \to C^(\D)$ for every $r> 0$ non-integer, as well as in $L^q(b\D) \to L^q(\D)$ for $1 \le q \le \infty$.

This is classical that properties of the Cauchy-Green operator $T$ in the Sobolev scale are similar to the regularity in the H\"older classes (see \cite{Ve}):
\begin{itemize}
\item[(i)] $T: W^{0,p}(\D) \to W^{1,p}(\D)$ is a bounded linear operator. 
\item[(ii)] $(Tf)_{\overline\zeta} = f$ i.e. $T$ solves the $\overline\partial$-equation in the unit disc (Sobolev's derivatives are used here). 
\item[(iii)] the function $Tf$ is holomorphic on $\C \setminus \overline\D$ and vanishes at infinity.
Furthemore, $Tf$ is $(1-2/p)$-H\"older continuous on $\C$  and the operator $T: W^{0,p}(\D) \to C^{1-2/p}(\C)$ is bounded.
\item[(iv)] $KTf(\zeta) = 0$ for $\zeta \in \D$.
\end{itemize}

Denote also by 
$\parallel f \parallel_\infty = \sup_\D \vert f \vert$
the usual $\sup$-norm on the space $L^\infty(\D)$ of complex functions bounded on $\D$.

Various versions of the following Lemma were used by several authors (see \cite{Fo, PiKh, NaRu}).

\begin{lemma}
\label{SchwarzLemma}
Let $f \in L^\infty(\D)$ and $f_{\overline\zeta} \in L^p(\D)$ for some $p > 2$. Then 
\begin{itemize}
\item[(a)]  $f$ admits a non-tangential limit at almost every point $\zeta \in b\D$.
\item[(b)] if $f$ admits a limit along a  curve in $\D$ approaching $b\D$ non-tangentially at a boundary point $e^{i\theta} \in b\D$, then $f$ admits a non-tangential limit at $e^{i\theta}$.
\item[(c)] for each positive $r < 1$ there exists a constant $C = C(r) > 0$ (independent of $f$) such that for every $\zeta_j \in r\D$, $j=1,2$ one has 
\begin{eqnarray}
\label{SchwarzIn}
\vert f(\zeta_1)  - f(\zeta_2)\vert \le C (\parallel f \parallel_\infty + \parallel f_{\overline\zeta} \parallel_{L^p(\D)} ) \vert \zeta_1 - \zeta_2\vert ^{1-2/p}
\end{eqnarray}
\end{itemize}
\end{lemma}
\proof The regularity property (i)  implies $Tf_{\overline\zeta} \in W^{1,p}(\D)$; in view of  (iii)   there exists $C_1 > 0$ independent of $f$ such that 
\begin{eqnarray}
\label{Sch1}
\vert Tf_{\overline\zeta}(\zeta_1)  - Tf_{\overline\zeta}(\zeta_2)\vert \le C_1 \parallel f_{\overline\zeta} \parallel_{L^p(\D)} \vert \zeta_1 - \zeta_2\vert ^{1-2/p}
\end{eqnarray}
The function $g = f - Tf_{\overline\zeta}$ is bounded on $\D$ and its generalized derivative vanishes: $g_{\overline\zeta} = 0$ on $\D$. Hence $g$ is holomorphic and (a), (b) follow respectively  from the classical Fatou and Lindel\"of theorems for holomorphic functions. Denote by $g^* \in L^\infty(b\D)$ and $f^*\in L^\infty(b\D)$ the non-tangential boundary value functions of $g$ and $f$ respectively. It is classical that $g$ satisfies the Cauchy formula on $\D$  that is $g(\zeta) =  Kg^*(\zeta) = Kf^*(\zeta) - KTf_{\overline\zeta}(\zeta) = 
Kf^*(\zeta)$ for each $\zeta \in \D$; here we have used the property (iv). Thus, the generalized Cauchy formula $f = K f^* + Tf_{\overline\zeta}$ holds on $\D$.
We have the estimate $\parallel Kf^* \parallel_\infty \le \parallel f \parallel_\infty$ and by the Cauchy estimates the holomorphic function $Kf^*$ is 1-Lipschitz on every $r\D$ with a Lipschitz constant $C(r)\parallel f \parallel_\infty$. In combination with (\ref{Sch1}) this proves (c).

\section{The Chirka-Lindel\"of principle for strictly pseudoconvex domains}

First we introduce an almost complex analog of an admissible approach which is classical in  the case of $\C^n$, see \cite{St,Ch}.

Let $\Omega$ be a smoothly bounded domain in an almost complex manifold $(M,J)$. Fix a hermitian metric on $M$ compatible with $J$; a choice of such metric will not affect our results since it changes only constant factors in estimates. We measure all distances and norms with  respect to the choosen metric.

Let $p \in b\Omega$ be a boundary point. {\it A non-tangential approach} to $b\Omega$ at $p$ can be defined as the limit along the sets
\begin{eqnarray}
\label{NonTan1}
C_\alpha(p) = \{ q \in \Omega: dist(q,p) < \alpha \delta_p(q) \}, \,\,\, \alpha > 1.
\end{eqnarray}
Here $\delta_p(q)$ denotes the minimum of distances from $q$ to the tangent plane $T_p(b\Omega)$ 
and to $b\Omega$.

We need to define a wider class of regions. {\it An admissible approach}  to $b\Omega$ at $p$ is defined as the limit along the sets

\begin{eqnarray}
\label{Ad1}
A_{\alpha}(p) = \{ q \in \Omega: d_p(q) < (1+\alpha)\delta_p(q), \, dist(p,q)^2 < \alpha \delta_p(q) \}, \,\, \alpha > 0.
\end{eqnarray}
Here $d_p(q)$ denotes the distance from $q$ to the holomorphic tangent space $H_p(b\Omega) = T_p(b\Omega) \cap J T_p(b\Omega)$. As in the classical case of $\C^n$, an admissible region approaches $b\Omega$ transversally in the normal direction and can be tangent in directions of the holomorphic tangent space.

\begin{definition}
A function $F:\Omega \to \C$ has an admissible limit $L$ at $p \in b\Omega$ if 
$\lim_{A_\alpha(p) \ni q} F(q) = L$ for all $\alpha > 0$.
\end{definition}

In what follows we denote by $b\D^+ = \{ e^{i\theta}: \theta \in [0,\pi] \}$ the upper-semi-circle.

\begin{definition}
Let $\Omega$ be a smoothly bounded domain in an almost complex manifold $(M,J)$ of complex dimension $n$. Assume that $f: \D \to \Omega$ is a $J$-complex disc of class $C^1(\overline\D)$ such that $f(b\D^+)$ is contained in $b\Omega$ and is transverse to $b\Omega$. Let also $\gamma: [0,1[ \to \D$ be a real curve of class $C^1([0,1])$, $\gamma(1) = i$ approaching $b\D$ non-tangentially at $i$.
Then the curve $\tau:= f \circ \gamma$ is called an admissible  $p$-curve, where $p = f(i) \in b\Omega$.
\end{definition}

\begin{definition}
A function $F$ defined on $\Omega$ has a limit $L \in \C$ along an admissible $p$-curve if there exists a $p$ curve $\tau$ such that $\lim_{t \to 1} (F \circ \tau)(t) = L$.
\end{definition}

As in the classical case, a smoothly bounded domain $\Omega$ in $(M,J)$ is    called {\it strictly pseudoconvex} if for every boundary point $p \in \Omega$ there exists a neighborhood $U$ of $p$ and a strictly $J$-plurisubharmonic function $\rho$ with non-vanishing gradient on $U$, such that $\Omega \cap U = \{ \rho < 0 \}$. Note that we do not need the existence of global defining strictly plurisubharmonic functions for $\Omega$ since all results are purely local.

Our first main result is the following analog of the Chirka - Lindel\"of principle \cite{Ch}.

\begin{thm}
\label{Thm1}
Let $\Omega$ be a smoothly bounded strictly pseudoconvex domain in an almost complex manifold $(M,J)$ of complex dimension $n$. Suppose that $F:\Omega \to \C$ is a bounded function of class $C^1(\Omega)$ and $\parallel\overline\partial_J F\parallel$ is  bounded on $\Omega$. If $F$ has a limit  along an admissible $p$-curve for some $p \in b\Omega$, then $F$ has the admissible limit  at $p$.
\end{thm}

Before proceed the proof we make some comments. Theorem \ref{Thm1} imposes two restrictions which are not present In the main result of E.Chirka \cite{Ch}. First, we assume that $\Omega$ is strictly pseudoconvex although \cite{Ch} deals with  much wider classes of domains. In particular, admissible regions in the sense of Chirka can be tangent to $b\Omega$ along complex submanifolds of 
higher dimension which can be contained in $b\Omega$. However, for a generic almost complex structure such complex submanifolds do not exist. Furthermore, any smooth boundary can be touched at a given point from inside by a strictly psedoconvex domain which allows to apply Theorem \ref{Thm1}. Second, in Chirka's theorem an admissible curve does not need to be contained in a complex disc. However, this assumption is sufficient for the most applications of Theorem \ref{Thm1}.

\bigskip

The local geometry of strictly pseudoconvex hypersurfaces in an almost complex manifold is similar to the case of $\C^n$ because of Lemma \ref{normalization}. However, there is some difference. In the case of $\C^n$ every strictly pseudoconvex hypersurface can be locally approximated by an osculating sphere or, equivalently, by the Siegel domain (the Heisenberg group). This is the standard and useful tool in local analysis on strictly pseudoconvex hypersurfaces. In the almost complex case the same remains true only in complex dimension 2. When the dimension is $> 2$, there exist an infinity local model almost complex structures which provides a local approximation; a choice of the model is determined by the first order jet of $J$ at a boundary point. All these model structures are strictly pseudoconvex and homogeneous. There are first appeared in \cite{GaSu} and later their geometry has been intensively studied.

We proceed the proof in several steps. First approximate a strictly pseudoconvex hypersurface by a suitable model structure. Then we establish Theorem for such a structure. The general result then follows by a perturbation argument.

\subsection{Local approximation by homogeneous models} To begin with, choose local coordinates near a boundary point according to Lemma \ref{normalization}. In these coordinates a local  defining function of $\Omega$ is strictly plurisubharmonic with respect to $J$ and $J_{st}$. It follows from the transformation rule (\ref{CompMat}) that the normalization conditions (\ref{norm}) are invariant with respect to usual (that is $J_{st}$) biholomorphic transformations.  This is well-known that using such 
(polynomial of degree at most two) transformations one can define $\Omega$ in a neighborhood of the origin by 
\begin{eqnarray}
\label{norm2}
\rho(z) = y_n + \vert 'z \vert^2 + o(\vert z \vert^2) < 0
\end{eqnarray}
where $'z = (z_1,....,z_{n-1})$. Notice that the normalization conditions (\ref{norm}) still hold in these coordinates.

Next for each $\lambda > 0$ we consider non-isotropic dilations 
\begin{eqnarray}
\label{dilations}
d_\lambda: ('z,z_n)  \mapsto (\lambda^{-1/2} {'z}, \lambda^{-1}z_n) = ('w,w_n).
\end{eqnarray}
The image $\Omega_\lambda:= d_\lambda(\Omega)$ is defined by 
$$\rho_\lambda:= \lambda^{-1}\rho(\lambda^{1/2}{'z},\lambda w_n) < 0$$
This is well-known (and is easy to check) that $\rho_\lambda$ converges uniformly on every compact subset of $\C^n$ to
the function  
\begin{eqnarray}
\label{model1}
\rho_0 = \Im w_n + \vert 'w \vert^2
\end{eqnarray}
as $\lambda \to 0$. 

Denote by $J_\lambda:= (d_\lambda)_*(J)$ the representation  of $J$ in the new coordinates; let also $A_\lambda(w)$ denotes the complex matrix of $J_\lambda$. Becuse of (\ref{norm}) the complex matrix of $J$ has the expansion $A(z) = L(\overline z) + O(\vert z \vert^2)$. Here $L$ is the linear part which depends only on $\overline z$.

 An elementary computation based on the transformation rule (\ref{CompMat}) shows that the functions $A_\lambda$ converge
uniformly on every compact subset of $\C^n$ as $\lambda \to 0$, to the matrix function 

\begin{eqnarray}
\label{model2}
A_0: ('w,w_n) \mapsto -\left(
\begin{array}{clll}
 0 & \dots & 0 & 0\\
 \dots & \dots & \dots & \dots\\
 0 & \dots & 0 & 0 \\
 l_1(\overline{'w})& \dots & l_{n-1}(\overline{'w}) & 0
\end{array} \right )
\end{eqnarray}
Here every $l_j$ in the last line is a complex linear function in $\overline{'w}$. Furthermore, every $l_j$ coincides with the restriction of the corresponding entry of the initial matrix function $L$ on the subspace $('w,0)$. Denote by $J_0$ the almost complex structure with the matrix $A_0$.

We call the domain (\ref{model1}) equipped with the almost complex structure   $J_0$ {\it a model structure}. The following properties of these structures are immediate:

\begin{itemize}
\item[(i)] every model structure is strictly pseudoconvex;
\item[(ii)] the non-isotropic dilations (\ref{dilations}) are  biholomorphic automorphisms of every model structure; in particular, each model structure is homogeneous.
\end{itemize}
Thus, model structures play the role of the Heisenberg group in the almost complex analysis on strictly pseudoconvex domains. 

In general, model structures are not integrable and there exists an infinity of biholomorphically non-equivalent model structures. The only exceptional case arises in dimension $2$ where all model structures are equivalent to $J_{st}$. Indeed, in this case the only non-zero entry of $A_0$ is $l_1(w_1) = a\overline{w}_1$ for some $a \in \C$. Using the transformation rule (\ref{CompMat}) one sees that the 
map $(w_1,w_2) \mapsto (w_1,w_2 + a\overline{w_1}^2/2)$ takes $A_0$ to $0$ that is $J_0$ becomes $J_{st}$. Obviously, after such a change of coordinates we obtain a domain biholomorphic to 
$\Omega_0$. Thus, in complex dimenson 2 the infinity of homogeneous models reduces to the usual ball with the standard complex structure.

\subsection{Case of model structures}  Here we establish the Chirka-Lindel\"of principle for the simplest case of model structures (\ref{model1}), (\ref{model2}). Note that for all these structures the complex normal line $\C \ni \zeta \mapsto ('0,\zeta)$ is $J_0$-complex. We will consider the case of an admissible $0$-curve which is contained in this line. We begin with the classical case of the $J_{st}$ which is the only one arising in dimension 2.

 Here the model structure is represented by the Siegel domain $\Omega_0$:

\begin{eqnarray}
\label{Heisn1}
\rho(z) = y_2 + \vert z_1 \vert^2 < 0
\end{eqnarray}
equipped with the standard complex structure $J_{st}$. Without loss of generality assume that the metric is the standard Euclidean (this choice affects only inessential constants in estimates). We use the notation $f(x) \sim g(x)$ for two functions $f(x),g(x)$ when there exists a constant $C > 0$ such that $C^{-1} g(x) \le f(x) \le C g(x)$. In what follows the value of constants $C$ can change from line to line.

We have  $T_0(b\Omega_0) = \{ y_2 = 0 \}$ and 
$H_0(b\Omega_0) = \{ z_2 = 0 \}$. Note that
 $$dist(z,b\Omega_0) \sim \vert \rho(z) \vert \le \vert y_2 \vert = dist(z,T_0(b\Omega_0)).$$ 
 Hence we can assume $\delta_0(z) = \vert\rho(z)\vert$. 
Since $dist(z,H_0(b\Omega_0)) = \vert z_2 \vert$, for each $\alpha > 0$ the admissible region $A_\alpha(0)$ is defined by the conditions

\begin{eqnarray}
\label{angle1}
\vert z_2 \vert < (1 + \alpha) \vert \rho(z) \vert
\end{eqnarray}
and
\begin{eqnarray}
\label{angle2}
\vert z \vert^2 < \alpha\vert \rho(z) \vert
\end{eqnarray}
The complex normal plane $(0,z_2)$ intersects $\Omega_0$ in the half-plane $\{(0,z_2): y_2 < 0 \}$ and the first inequality (\ref{angle1}) defines a non-tangential region there (which tends to this half-plane when $\alpha$ increases). 

Fix a point $(0,z_2^0)$ which satisfies (\ref{angle1}). Consider a complex line through the point $(0,z^0_2)$:

\begin{eqnarray}
\label{Heisn2}
f_a: \C \ni \zeta \mapsto (a\zeta, z^0_2)
\end{eqnarray}
and parallel to $H_0(b\Omega_0)$. A simple calculation shows that the second assumption  (\ref{angle2}) is equivalent to the fact that $f_a(\D) \subset A_\alpha(0)$ when
\begin{eqnarray}
\label{Heisn3}
\vert a \vert  \sim \vert y_2^0\vert^{1/2}
\end{eqnarray}
Clearly, this family of complex discs fills exactly the region $A_\alpha(0)$ when $(0,z_2^0)$ satisfies the first condition.

Next we consider the general case. Here the model structure is represented by the domain 
\begin{eqnarray}
\rho(z) = y_n + \vert 'z \vert^2 < 0
\end{eqnarray}
and the Cauchy-Riemann equations (\ref{holomorphy}) have the form

\begin{eqnarray}
\label{ModEq}
\left\{ \begin{array}{cccc}
& &(z_1)_{\overline\zeta} = 0, \\
& & ...\\
& & (z_{n-1})_{\overline\zeta} = 0, \\
& & (z_n)_{\overline\zeta} = \sum_{j=1}^{n-1} l_j(\overline{'z})(\overline{z}_j)_{\overline\zeta}
\end{array} \right.
\end{eqnarray}
The admissible regions are defined by the conditions

\begin{eqnarray}
\label{angle3}
\vert z_n \vert < (1 + \alpha) \vert \rho(z) \vert
\end{eqnarray}
and 
\begin{eqnarray}
\label{angle4}
 \vert z \vert^2 < \alpha\vert \rho(z) \vert
\end{eqnarray}
Consider any non-zero holomorphic tangent vector $v = (v_1,...,v_{n-1},0) \in H_0(b\Omega_0)$. Fix a point $('0,z_n^0)$ satisfying the assumptions (\ref{angle3}). Consider a pseudoholomorphic disc of the form

\begin{eqnarray}
\label{discs1}
\left\{ \begin{array}{cccc}
& & f(v,z_n^0): \D \ni \zeta \mapsto f(v,z_n^0)(\zeta), \\
& & 'z = v\zeta,\\
& & z_n = z_n^0 + l(v)\overline\zeta^2 + \overline{l(v)}\zeta^2.
\end{array} \right.
\end{eqnarray}
Here $l(v)$ is a suitable linear function in $v$ making $f$ a solution of the system (\ref{ModEq}). 

Note that the disc $f$ remains $J_0$-complex if we add  to its $z_n$-component any term holomorphic with respect to $J_{st}$. Our choice of this term does not affect the imaginary part of the map $z_n$ which remains equal to the constant $y_n^0$. Then as above, this disc is contained in $A_\alpha(0)$ when $\vert v \vert \sim \vert y_n^0\vert^{1/2}$ and the obtained family of discs fills the region $A_\alpha(0)$. After reparametrization we can assume that $\parallel v \parallel = 1$ and $f(v,z_n^0)$ is defined on the disc of radius $\sim \vert y_n^0 \vert^{1/2}$.

We sum up. Our construction of filling the admissible region $A_\alpha(0)$ consists of two steps. 

\begin{itemize} 
\item[(i)] First, 
we have the $J_0$-complex disc (in fact, this disc is parametrized by the half-plane $\Pi = \{ \Im \zeta < 0 \}$) of the form $f^0: \Pi \ni \zeta \mapsto ('0,\zeta)$. This disc is transverse to $b\Omega$; more precisely, it belongs to the complex normal line. 
\item[(ii)] Second, we have the family (\ref{discs1}) of 
$J_{0}$-complex discs $f(v,z_n^0)$ filling the region $A_\alpha(0)$. 
\item[(iii)] Every  map $f(v,z_n^0)$ is defined on the disc of radius $\sim \vert y_n^0\vert^{1/2}$, and is "parallel " (up to the second order terms) to $H_0(b\Omega)$; 
\item[(iv)] the centre $f(v,z_n^0)(0)$ of the disc is the point $('0,z_n^0)$ which lies in the non-tangential region in the transverse disc $f^0(\Pi)$;
\item[(v)] the unit tangent vector $v$ in the centre of the disc is parallel to $H_0(b\Omega)$.
\end{itemize}

Now we prove Theorem \ref{Thm1} for the case of a model structure. We suppose that $\Omega = \Omega_0$ with $J = J_{0}$. Furthermore, assume that the transverse disc $f$ coincides with the disc $f^0$ along the complex normal.  Let $F:\Omega_0 \to \C$ be a bounded $C^1$-function and 
$\parallel\overline\partial_{J_0} F \parallel$ also is bounded on $\Omega_0$.

 Then the restriction $F \circ f^0$ is a bounded function on $\Pi$ and $(F \circ f^0)_{\overline\zeta}$ is bounded as well. Furthermore, $F \circ f^0$ admits a limit $L$ along a non-tangential $0$-curve $\gamma$. By (b) Lemma \ref{SchwarzLemma} the function $F \circ f^0$ admits the limit $L$ along any non-tangential region in $\Pi$.
Let now $z \in A_\alpha(0)$. Then there exists (a unique) unit vector $v$ and a point $z_n^0$ in the non-tangential region on $f^0(\Pi)$ such that the disc $f(v,z^0_n)$  contains the point $z$ that is $z = f(v,z_n^0)(\zeta)$ for some $\zeta$ with $\vert \zeta \vert \le C \vert y_n^0 \vert^{1/2}$. Since also $f(v,z_n^0)(0) = z_n^0$, 
by (c) Lemma \ref{SchwarzLemma} we have the estimate (for any $0 < \tau < 1/2$): 
$$\vert F(z) - F('0,z_n^0) \vert = \vert (F \circ f(v,z_n^0))(\zeta) - (F \circ f(v,z_n^0))(0) \vert  \le C  \vert y_n^0 \vert^\tau$$
Since $F('0,z_n) \to L$ as $y_n^0 \to 0$, we conclude that $F(z) \to L$.

\subsection{Deformation argument and proof of Theorem \ref{Thm1}} In order to prove  Theorem \ref{Thm1} in general setting we construct a family of pseudoholomorphic discs with properties similar to the case of model structures. Locally every strictly pseudoconvex hypersurface is a small deformation of a model structure; a family of discs with required properties will in turn arise as a small deformation of discs constructed explicetly for the model case. We assume that we are in the hypothesis of Theorem; everything is local so $p = 0$.

First we consider the simplest case proving the existence of a limit along the non-tangential region 
$C_\alpha(0)$. Here the assumption of strict pseudoconvexity of $\Omega$ is irrelevant. As in the proof of Lemma \ref{lemma1}, consider isotropic dilations $h_\lambda$. One can assume that a local defining function of $\Omega$ has the form $y_n + o\vert z \vert $. We argue quite similarly to the case of model structures but using the isotropic dilations $h_\lambda$ instead on the non-isotropic dilations $d_\lambda$. Using the similar notations, we obtain a family of domains $\Omega_\lambda$ with 
almost complex structures $J_\lambda$ which for each $\lambda > 0$ are biholomorphic to $(\Omega,J)$.  When $\lambda \to 0$, these domains converge to the 
half-space $\Omega_0 = \{ y_n < 0 \}$ and $J_\lambda$ converge to $J_{st}$ in any $C^m$ norm.
Every non-tangential region for $\Omega_0$ is filled by a family of $J_{st}$-holomorphic discs as follows. First we have the transverse disc $f^0:\Pi \to \Omega$ defined as above. Second, we obtain the filling by the discs $\zeta \mapsto v\zeta + ('0,z_n^0)$ where $v = (v_1,...,v_{n-1},0)$ and $z_n^0$ lies in the non-tangent region in $\Pi$. Now, by the Nijenhuis-Woolf theorem, this construction is stable under small perturbation of the almost complex structure, that is, we have a family with similar properties for 
$(\Omega_\lambda,J_\lambda)$ for $\lambda > 0$ small enough. This gives the proof for $(\Omega_\lambda,J_\lambda)$ and so for $(\Omega,J)$ because they are biholomorphic.

\bigskip

Now we prove the existence of a limit in an admissible region $A_\alpha(0)$. There are at least two ways to proceed.  We will work with the structures $J_\lambda = (d_\lambda)_*(J)$ obtained from $J$ by non-isotropic dilations (\ref{dilations}); we also keep the notations $\Omega_\lambda$ introduced above. Let $f$ be a transverse $J$-complex disc given by hypothesis of Theorem \ref{Thm1}. It is more convenient to assume that $f$ is defined on the half-plane $\Pi = \{ \Im \zeta < 0 \}$, is of class $C^1$ up to the boundary and takes the segment $[-1,1]$ to $b\Omega$; after a holomorphic reparametrization we can assume $f(\zeta) = ('a,\zeta,\zeta) + o(\vert \zeta \vert)$. Then each disc $f_\lambda(\zeta):= (d_\lambda \circ f)(\lambda \zeta)$ is glued to $b\Omega_\lambda$ along $[-1/\lambda, 1/\lambda]$. The family $(f_\lambda)$ converges uniformly on compact subsets of $\overline \Pi$ to the $J_0$-complex  disc $f^0:\zeta \mapsto ('0,\zeta)$ glued to $\Omega_0$. This is exactly the one  considered in the model case. It suffices to prove the existence of a filling family of discs for $J_\lambda$ and $\Omega_\lambda$ when $\lambda$ is small enough. But this follows directly from the results on stability and deformation of pseudoholomorphic discs established in \cite{SuTu2}. Indeed, it is proved there (in a general setting) that a psedoholomorphic disc generates a family of nearby discs which are also pseudoholomorphic (for almost complex structures close enough to the initial one). The proof is somewhat inspired by the elementary argument of Nijenhuis-Woolf based on the implicit function theorem, but is more involved technically. It is based on an analysis of the Cauchy-Riemann equations (\ref{CRglobal}) linearized along a pseudoholomorphic disc; this provides the surjectivity of the linearized operator in a suitable scale (the H\"older or Sobolev one) necessary for the implicit function theorem. These stability result can be applied to the complex normal disc $f^0$ and to the family (\ref{discs1}) in complex tangent direction. This provides the existence of a filling family of disc for $(\Omega_\lambda),J_\lambda)$ and proves theorem.

The second approach is more elementary. First we observe that there are no problem with deformation of the "complex normal" disc $f^0$ since $J_0 \vert f^0 = J_{st}$. Hence the implicit function theorem can be applied similarly to the argument of Nijenhuis-Woolf. In order to construct a deformation of the family (\ref{discs1}) it suffices to apply the argument from Proposition 4.2 of \cite{DiSu}. It is  proved there that for each positive eigenvalue of the Levi form of the boundary there exists a pseudoholomorphic disc tangent to this direction and touching the boundary from outside at this point. By the Nijenhuis-Woolf theorem this disc generates a family of discs which enters to the domain in a prescribed holomorphic tangent direction at the centre. The stability of such a family follows again from the implicit function argument, see more details in \cite{DiSu}.

\section{The Fatou theorem}

As an application of the Chirka -Lindel\"of principle we establish a Fatou type results for $\overline\partial_J$-subsolutions. For holomorphic functions in $\C^n$ the first versions of the Fatou theorem are due to E.Stein \cite{St} and E.Chirka \cite{Ch}. Our approach is inspired by A.Sadullaev \cite{Sa1}.

We will deal with  some standard  classes of real submanifolds of an almost complex manifold. A submanifold $E$ of an almost complex $n$-dimensional $(M,J)$ is called 
{\it totally real} if at every point $p \in E$ the tangent space $T_pE$  does not contain non-trivial complex vectors that is $T_pE \cap JT_pE = \{ 0 \}$. This is well-known that the (real)  dimension of a totally real submanifold of $M$ is not bigger than $n$; we will consider in this paper only $n$-dimensional totally real submanifolds that is the case of maximal dimension. A real submanifold $N$ of $(M,J)$ is called {\it generic} if the complex span of $T_pN$ is equal to the whole  $T_pM$ for each point $p \in N$. A real $n$-dimensional submanifold of $(M,J)$ is generic if and only if it is totally real.

\bigskip

A totally real manifold $E$ can be defined as 
 \begin{eqnarray}
 \label{edge}
 E = \{ p \in M: \rho_j(p) = 0 \}
 \end{eqnarray}
 where $\rho_j:M   \to  \R$  are smooth functions with non-vanishing gradients.  The condition of total reality means that for every $p \in E$ the $J$-complex linear parts of the differentials $d\rho_j$ are (complex) linearly independent.

 A subdomain

\begin{eqnarray}
\label{wedge}
 W = \{ p \in M: \rho_j < 0, j = 1,...,n \}.
 \end{eqnarray}t
 is called {\it the wedge with the edge} $E$.

\bigskip

Our main result here is the following

\begin{thm}
\label{Thm2}
Let $E$ be a generic submanifold of the boundary $b\Omega$ of a smoothly bounded strictly pseudoconvex domain $\Omega$ in an almost complex manifold $(M,J)$ of complex dimension $n$. Suppose that $F:\Omega \to \C$ is a bounded function of class $C^1(\Omega)$ and 
$\parallel\overline\partial_J F\parallel$ is  bounded on $\Omega$. Then $F$ has an admissible limit  at almost every point of $E$.
\end{thm}
Note that the Hausdorff $n$-meausure on $E$ here is defined with respect to any metric on $M$; the condition to be a subset of measuro zero in $E$ is independent of such a choice.

In view of the following Lemma it suffices to consider the case where $E$ is totally real.

\begin{lemma}
\label{lemma2}
Let $N$ be a generic $(n+d)$-dimensional ($d > 0$) submanifold of an almost complex $n$-dimensional manifold $(M,J)$. Suppose that $K$ is a subset of $N$ of non-zero Hausdorff $(n+d)$-measure. Then 
there exists a (local) foliation of $N$ into a family $(E_s), s \in \R^{d}$ totally real $n$-dimensional submanifolds such that the intersection $K \cap E_s$ has a non-zero Hausdorff $n$-measure for 
each $s$ from some subset of non-zero Lebesgue measure in $\R^{d}$.
\end{lemma}
Here the Hausdorff measure is defined with respect to any Riemannian metric on $M$; the assumption that $K$ has a positive $n$-measure is independent on a choice of such metric.
\proof Let $p$ be a point of $M$ such that $K$ has a non-zero measure in each neighborhood of $p$. Choose local coordinates $z$ near  $p$ such that $p = 0$ and  $J(0) = J_{st}$.
After a $\C$-linear change of coordinates $N = \{ x_j + o(\vert z \vert) = 0, j = n-d+1,...,n\}$. After a local diffeomorphism with the identical linear part at $0$ we obtain that $N = \R^d(x_1,...,x_d) \times i\R^n(y)$. In the new coordinates the condition $J(0) = J_{s}$ still holds and every slice $E_s = \{ z \in N: x_1 = s_1,...,x_d= s_d\}$ is totally real. Now we conclude by the Fubini theorem.

\bigskip

Our proof of Theorem \ref{Thm2} uses the result of \cite{Su1}. There exist a wedge  $W \subset \Omega$ of the form (\ref{wedge})  with the edge $E$, and a family of $J$-complex discs $h_t: \D \to W$, of class $C^r(\overline\D)$ (with fixed $r > 1$) and smoothly depending on a parameter $t \in \R^{n-1}$, such that the following holds
\begin{itemize}
\item[(i)] $h_t(b\D^+) \subset E$ for every $t$; here $bD^+ = \{e^{i\theta}: \theta \in [0,\pi]\}$ is the upper semi-circle.
\item[(ii)] the curves $h_t(b\D^+)$ form a foliation of $E$. 
\item[(iii)] every $h_t(b\D^+)$ is transverse to $b\Omega$.
\item[(iv)] if $X_t \subset b\D^+$ is a subset of measure zero for every $t$, then $\cup_t h_t(X_t)$ ia a subset of measuro zero in $E$.
\end{itemize}
The function $F \circ h_t$ satisfies assumption of Lemma \ref{SchwarzLemma}. By (a) of this Lemma 
that for every $t$ the functon $F \circ h_t$  admits a non-tangential limit almost everywhere on $b\D^+$. Applying Theorem \ref{Thm1} we obtain that $F$ admits an admissible limt almost everywhere on $h_t(b\D^+)$. This implies theorem.

\bigskip

I conclude the paper by some remarks.

(1) Clearly, Theorems \ref{Thm1} and \ref{Thm2} are purely local that is they hold on an open strictly pseudoconvex piece of the boundary.

(2) As it was mentioned, the smoothness assumptions (on $M$, $J$, $b\Omega$) can be considerably weakened.

(3) It is quite probable that  more precise results can be obtained even using only  the methods of the present work, but I prefere to avoid technical complications. More advanced results will be considered in  forthcoming works.

(4) I hope that the methods of the present paper will allow to study boundary properties of other classes of functions 
on almost complex manifolds. One of the natural problems is to study boundary properties of almost holomorphic functions in the sense of S.Donaldson. To the best of my knowledge, the only result in this direction is due to M.Peyron \cite{Pe} who proved that a generic totally real manifold is a boundary uniqueness set for such functions.

\end{document}